\documentclass[11pt]{article}
\usepackage{amsmath,amsthm,amssymb,epsf}
\usepackage[dvips]{graphicx}

\theoremstyle{plain}

\theoremstyle{definition}

\theoremstyle{remark}

\DeclareMathOperator{\Ext}{Ext}

\newcommand{\QED}{\ifhmode\unskip\nobreak\fi\quad {\rm Q.E.D.}} 

\newcommand{\bQ}{\mathbb Q}

\newcommand\oM{\overline{M}}

\DeclareMathOperator{\vir}{vir}

\title{Open problems (for AGNES)}
\author{Rahul Pandharipande}
\date{April 2010}
\begin{document}
\maketitle

\setcounter{section}{-1}

Below are a few basic questions and speculations related to the moduli spaces
of curves, $K3$ surfaces, maps, and sheaves presented in 
the problem session of the AGNES conference in Amherst (April  2010).

\begin{center}
\begin{picture}(200,1)
\put(0,1){\line(1,0){200}}
\put(0,0){\line(1,0){200}}
\end{picture}
\end{center}

\vspace{12pt}
\noindent (i) {\em On the virtual class:}
\vspace{12pt}

Let $X$ be a nonsingular, projective variety over $\mathbb{C}$.
Let $\oM_g(X,\beta)$ be the moduli space of
 stable maps and let
$$\pi \colon \oM_g(X,\beta) \rightarrow \oM_g$$ be the forgetful morphism,
see \cite{fp} for background. The moduli space of stable
maps carries a virtual class 
$[\oM_g(X,\beta)]^{\vir}$
obtained from deformation theory 
\cite{Beh,BehFan,LiTian}. Tautological classes in the Chow and cohomology
rings of $\oM_g$ are defined efficiently in \cite{fabp}.
For a discussion of the properties (many conjectural) of tautological
classes, see also \cite{fab,panicm}.

\vspace{1cm}
\noindent
{\bf Q1}.  
Does $\pi_*[\oM_g(X,\beta)]^{\vir}\in H^*(\oM_g)$ 
 lie in the tautological ring in cohomology?
\vspace{1cm}

\noindent
{\bf Q2}.  
When does $\pi_*[\oM_g(X,\beta)]^{\vir}\in A^*(\oM_g)$ 
 lie in the tautological ring in Chow?
\vspace{1cm}

\noindent
I would guess the answer to ${\bf Q1}$ is yes. If $X$ is
a curve, an affirmative answer to ${\bf Q1}$ follows
from the results of \cite{fabp}.
We know $\pi_*[\oM_g(X,\beta)]^{\vir}$ does not always lie
in the tautological ring in Chow --- counterexamples can be found
when $X$ is a curve. A wild speculation, motivated
by the Bloch-Beilinson conjecture, is that the
answer to {\bf Q2} is yes when $X$ is defined over over $\overline{\bQ}$.

\vspace{18pt}
\noindent (ii) {\em On the Virasoro constraints:}
\vspace{12pt}

The spaces $\oM_g(X,\beta)$ determine the Gromov--Witten invariants of $X$. 
These are conjectured to satisfy the Virasoro constraints \cite{vir}.
Virasoro constraints are known to hold now in many, but not all, 
cases \cite{giv,opvir}.
A very interesting variety for which the Virasoro constraints are
unknown is the Enriques surface.

\vspace{1cm}
\noindent
{\bf Q3}. Prove the Virasoro constraints in case $X$ is an Enriques surface.
\vspace{1cm}

\noindent
A study of the Gromov-Witten theory of the Enriques surface,
closely related to modular forms, has been started in 
 \cite{mp2}. The Enriques surfaces
is perhaps the most basic variety where new techniques are
required to establish the Virasoro constraints.

\vspace{18pt}
\noindent (iii) {\em On the moduli of sheaves:}
\vspace{12pt}

Let $X$ be a nonsingular, projective 3-fold.
The Gromov-Witten theory of $X$, defined via
$\oM_g(X,\beta)$, is conjecturally \cite{mnop} equivalent to the
Donaldson-Thomas theory of $X$. The latter is defined
via the moduli of ideal sheaves of curves in $X$ \cite{DT,T}, or more recently,
in terms of the moduli spaces of stable pairs \cite{pt}.

\vspace{1cm}
\noindent
{\bf Q4}. Prove the GW/DT correspondence for $3$-folds.
\vspace{1cm}

\noindent
The toric cases of {\bf Q4} are known \cite{moop}. Algebraic cobordism
results \cite{levp} suggest the possibility of 
reducing to the toric case using degeneration methods.

\medskip

Donaldson--Thomas invariants are defined only in dimension $3$
because a virtual fundamental class for the moduli space of sheaves
is required. 
Deformations are given by $\Ext^1(E,E)$, obstructions by $\Ext^2(E,E)$, 
and to define the virtual fundamental class we need (roughly) the
vanishing
$$\Ext^i(E,E)=0 \ \  \text{for}\ \  i>2 \ .$$
On $3$-folds, the vanishing can often be obtained
 using Serre duality and stability. 
However, there are parallel examples of enumerative computations
in higher dimensions in Gromov-Witten theory \cite{klp,pz}. Moreover,
many aspects of Joyce's counting theory are valid in higher dimensions 
\cite{J}.

\vspace{1cm}
\noindent
{\bf Q5}. Define Donaldson--Thomas invariants in dimensions $>3$.
\vspace{.7cm}

\vspace{6pt}
\noindent (iv) {\em On the moduli of $K3$ surfaces:}
\vspace{12pt}

Let $M^{K3}_{2n}$ denote the moduli space of polarized $K3$ surfaces $(S,L)$ of 
degree $L^2=2n$. Little appears to be known about
the cycle theory of  $M^{K3}_{2n}$. 

\vspace{1cm}
\noindent
{\bf Q6}. What is the analogue of the tautological ring for $M^{K3}_{2n}$?
\vspace{1cm}

A natural guess for {\bf Q6} is the subring generated by the
classes of the 
Noether--Lefschetz loci. The Noether-Lefschetz loci parameterize
$K3$ surfaces with higher rank Picard lattices.

\vspace{1cm}
\noindent
{\bf Q7}. Do the Noether-Lefschetz divisors span 
 $\text{Pic}(M^{K3}_{2n}) \otimes_{\mathbb{Z}}{\mathbb{Q}}$?
\vspace{1cm}

Let $X$ be a compact 
Calabi-Yau 3-fold expressed as $K3$-fibration over $\mathbb{P}^1$,
$$\pi:X\rightarrow \mathbb{P}^1\ .$$
Given an ample line bundle $L$ on $X$, the family $\pi$
determines a morphism of the base $\mathbb{P}^1$ to the moduli
of polarized $K3$ surfaces.
Via \cite{mpnl}, the Gromov-Witten theory of $X$ in $\pi$-fiber classes is
calculated in terms of the Noether--Lefschetz numbers of $\pi$ and
the Katz-Klemm-Vafa \cite{kkv} conjecture concerning $\lambda_g$ integrals
in the reduced Gromov-Witten theory of a fixed $K3$ surface. The KKV
conjecture is proven for all classes in genus 0 in \cite{kmps} and
all genera in primitive classes in \cite{mpt}.

\vspace{1cm}
\noindent
{\bf Q8}. Prove the Katz-Klemm-Vafa conjecture for in all genera and
in all classes on $K3$ surfaces.
\vspace{1cm}

A solution to {\bf Q8} would provide a large class
of exact formulas for higher genus Gromov-Witten invariants of
 compact Calabi-Yau 3-folds. Unlike the local toric cases,
mathematical results for higher genus Gromov-Witten invariants 
have been difficult to obtain, see \cite{zing} for the genus 1
theory of the quintic 3-fold.

\vspace{1cm}
\noindent
{\bf Q9}. Find effective mathematical methods for calculating
the higher genus Gromov-Witten invariants of compact Calabi-Yau
3-folds.
\vspace{1cm}

Effective methods for the Enriques Calabi-Yau in genus $g\leq 2$
have been found in \cite{mp2}. Complete, but less effective, techniques for the
quintic are explained in \cite{mp1}. At present, 
the holomorphic anomaly equation
in topological string theory is more effective than the higher genus
mathematical methods.

\begin{center}
\begin{picture}(200,1)
\put(0,1){\line(1,0){200}}
\put(0,0){\line(1,0){200}}
\end{picture}
\end{center}

I was partially supported by the NSF through DMS-0500187.
Thanks are due to 
 P. Hacking and J. Tevelev for organizing the Algebraic
Geometry Northeastern Series (AGNES) conference. 
The present article was adapted from P. Hacking's notes
of my lecture.

$$ $$

\noindent Department of Mathematics \\
\noindent Princeton University

\end{document}